\documentclass[11pt]{article}
\usepackage{amsmath}
\usepackage{dsfont}
\usepackage{mathrsfs}
\usepackage{amsmath,amssymb}
\usepackage{amsfonts}
\usepackage{hyperref}
\usepackage{amsthm}
\usepackage{graphicx}
\usepackage{subfigure}
\usepackage{xcolor}
\usepackage{overpic}

\renewcommand{\qed}{\hfill\small{$\square$}\normalsize}

\hfuzz=\maxdimen
\theoremstyle{definition}
\newtheorem{lemma}{Lemma}[section]
\newtheorem{definition}[lemma]{Definition}
\newtheorem{proposition}[lemma]{Proposition}
\newtheorem{theorem}[lemma]{Theorem}

\newtheorem{remark}{Remark}

\numberwithin{equation}{section}
\renewcommand{\proof}{\textbf{Proof. }}
\renewcommand{\qed}{\hfill\small{$\square$}\normalsize}

\DeclareFixedFont{\Acknowledgment}{OT1}{cmr}{bx}{n}{14pt}
\begin{document}

\title{\bf On the global rigidity of sphere packings on 3-dimensional manifolds}
\author{Xu Xu}
\date{}
\maketitle

\begin{abstract}
In this paper, we prove the global rigidity of sphere packings
on 3-dimensional manifolds.
This is a 3-dimensional analogue of the rigidity theorem of Andreev-Thurston and was conjectured by Cooper and Rivin in \cite{CR}.
We also prove a global rigidity result using a combinatorial scalar curvature introduced by
Ge and the author in \cite{GX4}.
\end{abstract}



\section{Introduction}\label{section 1}
In his investigation of hyperbolic metrics on 3-manifolds,
Thurston (\cite{T1}, Chapter 13) introduced the circle packing with prescribed intersection angles
and proved the Andreev-Thurston Theorem, which consists of two parts.
The first part is on the existence of circle packing for a given triangulation. The second part is on the rigidity of circle packings,
which states that a circle packing is uniquely determined by its discrete Gauss curvature (up to scaling for the Euclidean background geometry).
For a proof of Andreev-Thurston Theorem, see \cite{CL1, DV, H, MR, S, T1}.

To study the $3$-dimensional analogy of the circle packing on surfaces,
Cooper and Rivin \cite{CR} introduced the sphere packing on 3-dimensional manifolds.
Suppose $M$ is a 3-dimensional closed manifold with a triangulation $\mathcal{T}=\{V,E,F,T\}$,
where the symbols $V,E,F,T$
represent the sets of vertices, edges, faces and tetrahedra, respectively.
\begin{definition}[\cite{CR}]\label{definition of sphere packing metric}
A Euclidean (hyperbolic respectively) sphere packing metric on $(M, \mathcal{T})$ is a map $r:V\rightarrow (0,+\infty)$
such that (1) the length of an edge $\{ij\}\in E$ with vertices $i, j$ is $l_{ij}=r_{i}+r_{j}$
and (2) for each tetrahedron $\{i,j,k,l\}\in T$,
the lengths $l_{ij},l_{ik},l_{il},l_{jk},l_{jl},l_{kl}$ form the edge lengths of a Euclidean (hyperbolic respectively) tetrahedron.
\end{definition}

The condition (2) is called the nondegenerate condition, which makes the space of sphere packing metrics to be a proper
open subset of $\mathbb{R}^{|V|}_{>0}$.
The space of sphere packing metrics will be denoted by $\Omega$ in the paper.

To study sphere packing metrics,
Cooper and Rivin \cite{CR} introduced the combinatorial scalar curvature $K: V\rightarrow \mathbb{R}$,
which is defined as angle deficit of solid angles at a vertex $i$
\begin{equation}\label{CR curvature}
K_{i}= 4\pi-\sum_{\{i,j,k,l\}\in T}\alpha_{ijkl},
\end{equation}
where $\alpha_{ijkl}$ is the solid angle at the vertex $i$ of the tetrahedron $\{i,j,k,l\}\in T$
and the summation is taken over all tetrahedra with $i$ as a vertex.

The sphere packing metrics have the following local rigidity with respect to the combinatorial scalar curvature $K$.

\begin{theorem}[\cite{CR, G1, G2,R}]
Suppose $(M, \mathcal{T})$ is a closed $3$-dimensional triangulated manifold.
Then
a Euclidean or hyperbolic sphere packing metric on $(M, \mathcal{T})$ is
locally determined by its combinatorial scalar curvature $K$ (up to scaling for the Euclidean background geometry).
\end{theorem}

The global rigidity of sphere packing metrics on 3-dimensional triangulated manifolds was conjectured by Cooper and Rivin in \cite{CR}.
In this paper, we solve this conjecture and prove the following result.

\begin{theorem}\label{main theorem rigidity}
Suppose $(M, \mathcal{T})$ is a closed triangulated $3$-manifold.
\begin{description}
  \item[(1)] A Euclidean sphere packing metric on $(M, \mathcal{T})$ is determined by its combinatorial scalar curvature $K: V\rightarrow \mathbb{R}$ up to scaling.
  \item[(2)] A hyperbolic sphere packing metric on $(M, \mathcal{T})$ is determined by its combinatorial scalar curvature $K: V\rightarrow \mathbb{R}$.
\end{description}
\end{theorem}

Although the combinatorial curvature $K$ is a good candidate for the 3-dimensional combinatorial scalar curvature,
it has two disadvantages comparing to the smooth scalar curvature on Riemannian manifolds.
The first is that it is scaling invariant with respect to the Euclidean sphere packing metrics, i.e. $K(\lambda r)=K(r)$ for $\lambda>0$;
The second is that $K_i$ tends to zero as the triangulation of the manifold is finer and finer.
Motivated by the observations, Ge and the author \cite{GX4} introduced a new combinatorial scalar curvature defined as
$R_i=\frac{K_i}{r_i^2}$ for 3-dimensional manifolds with Euclidean background geometry, which overcomes the two disadvantages
if we take $g_i=r_i^2$ as an analogue of the Riemannian metric tensor for the Euclidean background geometry.
This definition can be modified to fit the case of hyperbolic background geometry.
We further generalized this definition of combinatorial scalar curvature to the following combinatorial $\alpha$-curvature.
\begin{definition}[\cite{GX4, GX5}]\label{definition of alpha curvature}
Suppose $(M, \mathcal{T})$ is a closed triangulated $3$-manifold with a sphere packing metric $r: V\rightarrow (0, +\infty)$ and $\alpha\in \mathbb{R}$.
Combinatorial $\alpha$-curvature at a vertex $i\in V$ is defined to be
\begin{equation}\label{definition of R-curvature}
R_{\alpha,i}=\frac{K_i}{s_i^{\alpha}},
\end{equation}
where $s_i=r_i$ for the Euclidean sphere packing metrics and $s_i=\tanh \frac{r_i}{2}$ for the hyperbolic sphere packing metrics.
\end{definition}
When $\alpha=0$, the $0$-curvature $R_0$ is the combinatorial scalar curvature $K$.
When $\alpha=-1$, $R_{-1,i}=K_ir_i$ is closely related to the discrete curvature given by Regge \cite{Re}
as described by Glickenstein in Section 5.1 of \cite{G4}.
Regge's formulation is shown to converge to scalar curvature measure $RdV$ by
Cheeger-M\"{u}ller-Schrader \cite{CMS}, which indicates that $K_ir_i$ is an analogue of $RdV$.
Combinatorial $\alpha$-curvatures on triangulated surfaces were studied in \cite{GJ3, GX2, GX4, GX5, GX3, GX6, X}.
Using the combinatorial $\alpha$-curvature, we prove the following global rigidity on 3-manifolds.

\begin{theorem}\label{main theorem global rigidity for alpha curvature}
Suppose $(M, \mathcal{T})$ is a closed triangulated $3$-manifold and $\overline{R}$
is a given function defined on the vertices  of  $(M, \mathcal{T})$.
\begin{description}
  \item[(1)]  In the case of Euclidean background geometry,
  \begin{description}
    \item[(a)] if $\alpha\overline{R}\equiv0$, there exists at most one Euclidean sphere packing metric  in $\Omega$
  with combinatorial $\alpha$-curvature equal to $\overline{R}$ up to scaling.
    \item[(b)] if $\alpha\overline{R}\leq0$ and $\alpha\overline{R}\not\equiv0$,
  there exists at most one Euclidean sphere packing metric in $\Omega$ with combinatorial $\alpha$-curvature equal to $\overline{R}$.
  \end{description}

  \item[(2)] In the case of hyperbolic background geometry, if $\alpha\overline{R}\leq 0$,
  there exists at most one hyperbolic sphere packing metric in $\Omega$
  with combinatorial $\alpha$-curvature equal to $\overline{R}$.
\end{description}
\end{theorem}

When $\alpha=0$, Theorem \ref{main theorem global rigidity for alpha curvature}
is reduced to Theorem \ref{main theorem rigidity}.
When $\alpha=2$, the local rigidity of Euclidean sphere packing metrics with nonpositive
constant $2$-curvature was proven in \cite{GX4}.
For $\alpha\in\mathbb{R}$, the local rigidity of Euclidean sphere packing metrics with constant
combinatorial $\alpha$-curvature on 3-dimensional triangulated manifolds
was proven in \cite{GX5}. Results similar to Theorem \ref{main theorem global rigidity for alpha curvature}
were proven for Thurston's circle packing metrics on surfaces in \cite{GX4, GX3}
and for inversive distance circle packing metrics on surfaces in \cite{GJ1, GJ2, GJ3, GX6, X}.

Glickenstein \cite{G1} introduced a combinatorial Yamabe flow to study the
constant curvature problem of $K$. He found that the combinatorial scalar curvature $K$
evolves according to a heat type equation along his flow and showed that the solution converges to a constant curvature
metric under some nonsingular conditions.
Glickenstein \cite{G2} further derived a maximal principle for the curvature along the combinatorial Yamabe flow under certain assumptions on the triangulation.
Ge and the author \cite{GX2, GX4, GX5} generalized Cooper and Rivin's definition of combinatorial scalar curvature
and introduced a combinatorial Yamabe flow to deform the sphere packing metrics, aiming at
finding the corresponding constant curvature sphere packing metrics on 3-dimensional triangulated manifolds.
Ge and Ma \cite{GM} studied the deformation of combinatorial $\alpha$-curvature on $3$-dimensional triangulated manifolds
using a modified combinatorial Yamabe flow.

The paper is organized as follows.
In Section \ref{Section 2}, We give a description of the admissible space of sphere packing metrics for a single tetrahedron.
In Section \ref{Section 3}, we recall Cooper and Rivin's action functional and extend it to be a convex functional.
In Section \ref{Section 4}, We prove Theorem \ref{main theorem rigidity} and Theorem \ref{main theorem global rigidity for alpha curvature}.

\section{Admissible space of sphere packing metrics for a single tetrahedron}\label{Section 2}
Suppose $M$ is a 3-dimensional connected closed manifold with a triangulation $\mathcal{T}=\{V,E,F,T\}$.
We consider sphere packing metrics as points in $\mathbb{R}^N_{>0}$,
where $N=|V|$ denotes the number of vertices.
And we use $\mathbb{R}^V$ to denote the set of real functions defined on the set of vertices $V$.

Suppose $r$ is a Euclidean sphere packing metric on $(M, \mathcal{T})$.
For any edge $\{ij\}\in E$, let $l_{ij}=r_i+r_j$
and for a tetrahedron $\{ijkl\}\in T$, $l_{ij}, l_{ik}, l_{il}, l_{jk}, l_{jl}, l_{kl}$ can be realized as edge lengths of
a Euclidean tetrahedron.
Gluing all of these Euclidean tetrahedra in $T$ along the faces isometrically
produces a piecewise linear metric on the triangulated manifold $(M, \mathcal{T})$.
On this manifold, drawing a sphere $S_i$ centered at vertex $i$ of radius $r_i$ for each vertex $i\in V$,
we obtain a Euclidean sphere packing. A hyperbolic sphere packing can be constructed similarly.

As $l_{ij}=r_i+r_j$, it is straightforward to show that the triangle inequalities for $l_{ij}, l_{ik}, l_{jk}$ hold on the face $\{ijk\}\in F$.
However, triangle inequalities on the faces are not enough for
$l_{ij}, l_{ik}, l_{il}, l_{jk}, l_{jl}, l_{kl}$ to determine a Euclidean or hyperbolic tetrahedron.
There are nondegenerate conditions.
It is found \cite{CR, G1}
that Descartes circle theorem, also called Soddy-Gossett theorem,
can be used to describe the degenerate case. We state a version obtained in \cite{M}.

An oriented circle is a circle together with an assigned direction of unit normal. The interior of an oriented circle
is its interior for an inward pointing normal and its exterior for an outward pointing normal.

\begin{definition}[\cite{LMW}]
A Euclidean (hyperbolic respectively) oriented Descartes configuration consists of
$4$ mutually tangent oriented circles in the Euclidean (hyperbolic respectively) plane
such that all pairs of tangent circles have distinct points of tangency and the interiors of all four oriented circles
are disjoint.
\end{definition}

\begin{figure}[!htb]
\centering
  \includegraphics[height=0.43\textwidth,width=1\textwidth]{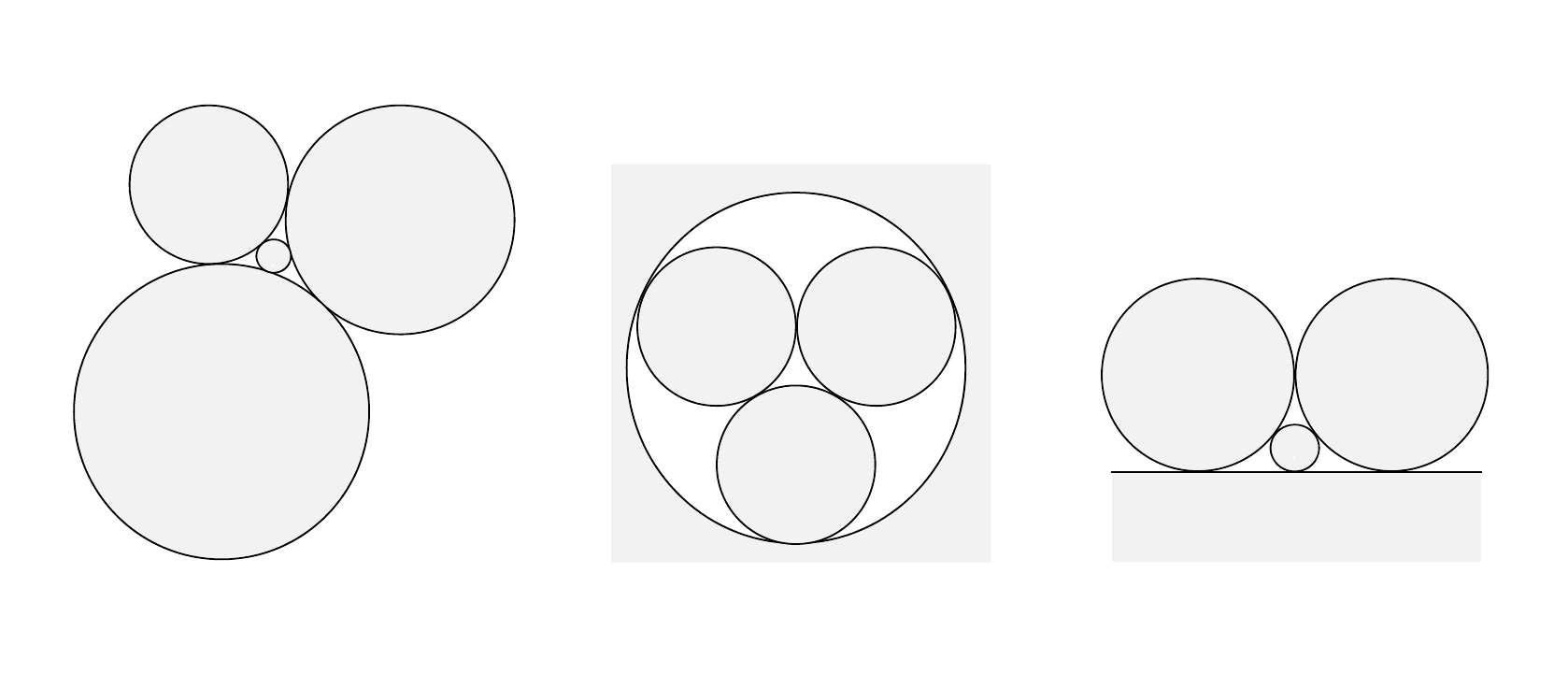}
  \caption{Descartes configurations}
  \label{Descartes configurations}
\end{figure}

Several Euclidean oriented Descartes configurations are shown in Figure \ref{Descartes configurations},
where the shadow denotes the interior of a circle.

\begin{remark}
We allow the Euclidean oriented Descartes configuration to include the straight lines and
the hyperbolic oriented Descartes configuration to include the horocycles.
\end{remark}

\begin{theorem}[Descartes Circle Theorem]\label{Descartes circle theorem}
~

\begin{description}
  \item[(1)] Given a Euclidean Descartes configuration $C_i, i=1,2,3,4$ such that
$C_i$ has radius $r_i$. Then
\begin{equation*}
\left(\sum_{i=1}^4k_i\right)^2-2\sum_{i=1}^4k_i^2=0,
\end{equation*}
where $k_i=\frac{1}{r_i}$, if $C_i$ is assigned an inward pointing normal, otherwise $k_i=-\frac{1}{r_i}$.
  \item[(2)] Given a hyperbolic Descartes configuration
  $C_i, i=1,2,3,4$ such that $S_i$ has radius $r_i$. Then
\begin{equation*}\label{hyperbolic quadratic}
\begin{aligned}
\left(\sum_{i=1}^4k_i\right)^2-2\sum_{i=1}^4k_i^2+4=0,
\end{aligned}
\end{equation*}
where $k_i=\coth r_i$, if $C_i$ is assigned an inward pointing normal, otherwise $k_i=-\coth r_i$.
\end{description}
\end{theorem}

There is also a version of Descartes circle theorem for the spherical background geometry.
See \cite{LMW, M} for Descartes circle theorem with different background geometries.
In this paper, we concentrate on the Euclidean and hyperbolic cases.

For the Euclidean background geometry,
Glickenstein \cite{G1} observed the admissible space of Euclidean sphere packing metrics for a tetrahedron $\{ijkl\}\in T$ to be a Euclidean tetrahedron is
\begin{equation*}
\Omega^{\mathbb{E}}_{ijkl}=\{(r_i, r_j, r_k, r_l)\in \mathbb{R}^4_{>0}|Q^{\mathbb{E}}_{ijkl}>0\},
\end{equation*}
where
\begin{equation*}
Q^{\mathbb{E}}_{ijkl}=\left(\frac{1}{r_{i}}+\frac{1}{r_{j}}+\frac{1}{r_{k}}+\frac{1}{r_{l}}\right)^2-
2\left(\frac{1}{r_{i}^2}+\frac{1}{r_{j}^2}+\frac{1}{r_{k}^2}+\frac{1}{r_{l}^2}\right).
\end{equation*}
For the hyperbolic background geometry, we need the following result.
\begin{proposition}[\cite{TW}, Proposition 2.4.1]\label{nondegenerate condition for hyperbolic}
A non-degenerate hyperbolic tetrahedron with edge lengths $l_{ij}, l_{ik}, l_{il}, l_{jk}, l_{jl}, l_{kl}$ exists
if and only if all principal minors of
\begin{equation*}
\begin{aligned}
\left|
                   \begin{array}{ccccc}
                     1 & \cosh l_{ij} & \cosh l_{ik} & \cosh l_{il}\\
                    \cosh l_{ij} & 1 & \cosh l_{jk} & \cosh l_{jl}\\
                     \cosh l_{ik} & \cosh l_{jk} & 1 & \cosh l_{kl}\\
                     \cosh l_{il} & \cosh l_{jl} & \cosh l_{kl} & 1 \\
                   \end{array}
                   \right|
\end{aligned}
\end{equation*}
are negative.
\end{proposition}

Applying Proposition \ref{nondegenerate condition for hyperbolic} to hyperbolic sphere packing metrics, we have
the admissible space of hyperbolic sphere packing metrics
for a tetrahedron $\{ijkl\}\in T$ to be a non-degenerate hyperbolic tetrahedron is
\begin{equation*}
\Omega^{\mathbb{H}}_{ijkl}=\{(r_i, r_j, r_k, r_l)\in \mathbb{R}^4_{>0}|Q^{\mathbb{H}}_{ijkl}>0\},
\end{equation*}
where
\begin{equation*}
\begin{aligned}
Q^{\mathbb{H}}_{ijkl}=&\left(\coth r_i+\coth r_j+\coth r_k+\coth r_l\right)^2\\
&-2\left(\coth^2 r_i+\coth^2 r_j+\coth^2 r_k+\coth^2 r_l\right)+4.
\end{aligned}
\end{equation*}

Cooper and Rivin \cite{CR} called the tetrahedra produced by sphere packing conformal and
proved that a tetrahedron is a Euclidean conformal tetrahedron if and only if there exists a unique sphere
tangent to all of the edges of the tetrahedron.
Moreover, the point of tangency with the
edge $\{ij\}$ is of distance $r_i$ to $i$-th vertex.
They proved the following lemma on the admissible space of sphere packing metrics for a single tetrahedron.
\begin{lemma}[\cite{CR}]\label{simply connectness of ijkl}
For a Euclidean or hyperbolic tetrahedron $\{ijkl\}\in T$,
the admissible spaces $\Omega^{\mathbb{E}}_{ijkl}$ and $\Omega^{\mathbb{H}}_{ijkl}$ are simply connected open subsets of $\mathbb{R}^4_{>0}$.
\end{lemma}
They further pointed out that $\Omega^{\mathbb{E}}_{ijkl}$ is not convex.
For a triangulated $3$-manifold $(M, \mathcal{T})$,
the admissible spaces
\begin{equation*}
\Omega^{\mathbb{E}}=\{r\in \mathbb{R}^N_{>0}|Q_{ijkl}^{\mathbb{E}}>0, \forall \{ijkl\}\in T\}
\end{equation*}
and
\begin{equation*}
\Omega^{\mathbb{H}}=\{r\in \mathbb{R}^N_{>0}|Q^{\mathbb{H}}_{ijkl}>0, \forall \{ijkl\}\in T\}
\end{equation*}
are open subsets of $\mathbb{R}^N_{>0}$.

We need a good description of the admissible spaces $\Omega^{\mathbb{E}}_{ijkl}$ and $\Omega^{\mathbb{H}}_{ijkl}$ for a single tetrahedron $\{ijkl\}\in T$.
If the radii $r_j, r_k, r_l$ of the spheres $S_j, S_k, S_l$ are fixed, Cooper and Rivin \cite{CR} observed that
degeneracy occurs when $r_i$ is large enough so that the sphere $S_i$ is large enough to be tangent to the other
three spheres, yet small enough that its center $i$ lies in the plane defined by $j, k$ and $l$.
This defines a degenerate set $V_i$ of the sphere packing metrics.
The degenerate sets $V_j, V_k$ and $V_l$ can be defined similarly.

We have the following result on the structure of the sets $V_i$, $V_j$, $V_k$, $V_l$ and $\Omega_{ijkl}$.
Here and in the rest of the paper,
$\Omega_{ijkl}$ denotes $\Omega^{\mathbb{E}}_{ijkl}$ or $\Omega^{\mathbb{H}}_{ijkl}$ according to the background geometry
and $\overline{\Omega}_{ijkl}$ denotes the closure of $\Omega_{ijkl}$ in $\mathbb{R}^4_{>0}$.

\begin{theorem}\label{structure of admissible r}
Connected components of $\mathbb{R}^4_{>0}-\Omega_{ijkl}$ are $V_i$, $V_j$, $V_k$ and $V_l$.
Furthermore, the intersection of $\overline{\Omega}_{ijkl}$ with any of $V_i$, $V_j$, $V_k$, $V_l$ is a connected component of
$\overline{\Omega}_{ijkl}-\Omega_{ijkl}$, which is
a graph of a continuous and piecewise analytic function defined on $\mathbb{R}^3_{>0}$.
\end{theorem}
\proof
We prove the Euclidean case in details. The hyperbolic case is similar and will be omitted.
By symmetry, it suffices to consider $V_i$.
It is observed \cite{G2} that
\begin{equation*}
\begin{aligned}
Q^{\mathbb{E}}_{ijkl}
=&\frac{1}{r_i}(\frac{1}{r_j}+\frac{1}{r_k}+\frac{1}{r_l}-\frac{1}{r_i})+\frac{1}{r_j}(\frac{1}{r_i}+\frac{1}{r_k}+\frac{1}{r_l}-\frac{1}{r_j})\\
  &+\frac{1}{r_k}(\frac{1}{r_i}+\frac{1}{r_j}+\frac{1}{r_l}-\frac{1}{r_k})+\frac{1}{r_l}(\frac{1}{r_i}+\frac{1}{r_j}+\frac{1}{r_k}-\frac{1}{r_l}).
\end{aligned}
\end{equation*}
If $r_i=\min\{r_i, r_j, r_k, r_l\}$ and $Q^{\mathbb{E}}_{ijkl}=0$, then
$\frac{1}{r_j}+\frac{1}{r_k}+\frac{1}{r_l}-\frac{1}{r_i}<0$ and
\begin{equation}\label{derivative of Q}
\frac{\partial Q^{\mathbb{E}}_{ijkl}}{\partial r_i}=-\frac{2}{r_i^2}(\frac{1}{r_j}+\frac{1}{r_k}+\frac{1}{r_l}-\frac{1}{r_i})>0.
\end{equation}
This implies that if $Q^{\mathbb{E}}_{ijkl}=0$, we can always increase $r_i$ to make the tetrahedron nondegenerate.
So we just need to analyze the critical degenerate case of $V_i$,
where $S_i$ is externally tangent to the other three spheres $S_j, S_k, S_l$
and the center $i$ lies in the plane defined by $j, k$ and $l$.
Glickenstein further proved the following result.

\begin{proposition}[\cite{G2}, Proposition 6]\label{degenerate lemma}
If $Q^{\mathbb{E}}_{ijkl}\rightarrow 0$ such that none of $r_i, r_j, r_k, r_l$ tend to $0$,
then one solid angle tends to $2\pi$ and
the others tend to 0.
Furthermore, if $r_i$ is the minimum of $r_i, r_j, r_k, r_l$, then the dihedral angles $\beta_{ijkl}, \beta_{ikjl}, \beta_{iljk}$
tend to $\pi$, the dihedral angles $\beta_{jkil}, \beta_{jlik}, \beta_{klij}$ tend to $0$ and the solid angle $\alpha_{ijkl}$ tends
to $2\pi$, where $\beta_{ijkl}$ is the dihedral angle along the edge $\{ij\}$.
\end{proposition}

\begin{figure}[!htb]
\centering
  \includegraphics[height=0.47\textwidth,width=1\textwidth]{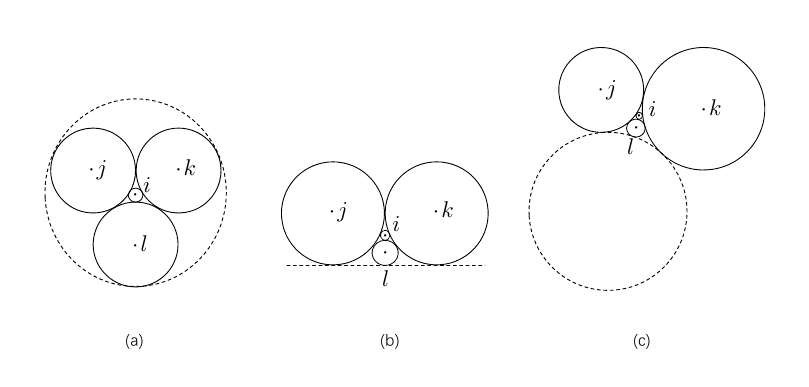}
  \caption{circle configurations (the circles with dotted line correspond to the roots of (\ref{quadratic equation in r_i}) rejected)}
  \label{circle}
\end{figure}

Proposition \ref{degenerate lemma} implies that $V_i, V_j, V_k, V_l$ are the only degenerations that can occur.
Furthermore, if $Q^{\mathbb{E}}_{ijkl}=0$ and $r_i<\min\{r_j, r_k, r_l\}$,
the center $i$ lies in the interior of the Euclidean triangle $\triangle jkl$,
which determines an oriented Descartes configuration
consisting of four externally tangent circles $C_i, C_j, C_k, C_l$ in the plane defined by $j, k$ and $l$.
See Figure \ref{circle}.

\begin{lemma}
Suppose $C_i, C_j, C_k, C_l$ are four oriented circles with finite radii and inward pointing normal in the Euclidean plane.
If $C_i, C_j, C_k, C_l$ form an oriented Descartes configuration and $r_i<\min\{r_j, r_k, r_l\}$, then
\begin{equation}\label{equation of r_i}
\begin{aligned}
r_i=f(r_j,r_k,r_l):=\left\{
      \begin{array}{ll}
        \frac{-B+\sqrt{B^2-4AC}}{2A}, & \hbox{$(r_j, r_k, r_l)\in \Omega_{jkl}$;} \\
        -\frac{C}{B}, & \hbox{$(r_j, r_k, r_l)\in \overline{\Omega}_{jkl}\setminus\Omega_{jkl}$;} \\
        \frac{-B+\sqrt{B^2-4AC}}{2A}, & \hbox{$(r_j, r_k, r_l)\in \mathbb{R}^3_{>0}\setminus \overline{\Omega_{jkl}}$;}
      \end{array}
    \right.
\end{aligned}
\end{equation}
where
\begin{equation*}
\begin{aligned}
A=&2r_jr_kr^2_l+2r_jr_lr^2_k+2r_kr_lr^2_j-r_k^2r_l^2-r_j^2r_l^2-r_j^2r_k^2,\\
B=&2 r_jr_kr_l(r_kr_l+r_jr_l+r_jr_k),\\
C=&-r_j^2r_k^2r_l^2,\\
\Omega_{jkl}=&\{(r_j,r_k,r_l)\in \mathbb{R}^3_{>0}|A>0\}.
\end{aligned}
\end{equation*}
\end{lemma}
\proof By Descartes circle theorem \ref{Descartes circle theorem}, we have
\begin{equation*}
\begin{aligned}
Q^{\mathbb{E}}_{ijkl}=\left(\frac{1}{r_{i}}+\frac{1}{r_{j}}+\frac{1}{r_{k}}+\frac{1}{r_{l}}\right)^2 -2\left(\frac{1}{r_{i}^2}+\frac{1}{r_{j}^2}+\frac{1}{r_{k}^2}+\frac{1}{r_{l}^2}\right)=0,
\end{aligned}
\end{equation*}
which is equivalent to the quadratic equation in $r_i$
\begin{equation}\label{quadratic equation in r_i}
\begin{aligned}
Ar_i^2+Br_i+C=0.
\end{aligned}
\end{equation}
For the quadratic equation (\ref{quadratic equation in r_i}) in $r_i$, the discriminant is
\begin{equation*}
\begin{aligned}
\Delta=&B^2-4AC=16r_j^3r_k^3r_l^3(r_j+r_k+r_l),
\end{aligned}
\end{equation*}
which is always positive for $(r_j, r_k, r_l)\in \mathbb{R}^3_{>0}$.
Note that
\begin{equation*}
\begin{aligned}
A=&2r_jr_kr^2_l+2r_jr_lr^2_k+2r_kr_lr^2_j-r_k^2r_l^2-r_j^2r_l^2-r_j^2r_k^2\\
=&(\sqrt{r_jr_k}+\sqrt{r_jr_l}+\sqrt{r_kr_l})(\sqrt{r_jr_k}+\sqrt{r_jr_l}-\sqrt{r_kr_l})\\
&(\sqrt{r_jr_k}-\sqrt{r_jr_l}+\sqrt{r_kr_l})(-\sqrt{r_jr_k}+\sqrt{r_jr_l}+\sqrt{r_kr_l}).
\end{aligned}
\end{equation*}

\begin{description}
  \item[(1)] In the case that $A>0$, we have $\sqrt{r_jr_k}, \sqrt{r_jr_l}, \sqrt{r_kr_l}$ satisfy the triangle inequalities. As $B>0$ and $C<0$,
the quadratic equation (\ref{quadratic equation in r_i}) has two roots with different signs and negative sum.
Note that $r_i>0$, by Theorem \ref{Descartes circle theorem}, we have
$$r_i=\frac{-B+\sqrt{\Delta}}{2A}>0.$$
The negative root corresponds to the case that the mutually externally tangent circles $C_j, C_k, C_l$ are internally tangent to a circle $C'_i$.
See (a) in Figure \ref{circle}.

  \item[(2)] In the case that $A=0$, i.e. $\sqrt{r_jr_k}+\sqrt{r_jr_l}=\sqrt{r_kr_l}$ or $\sqrt{r_jr_k}+\sqrt{r_kr_l}=\sqrt{r_jr_l}$
or $\sqrt{r_kr_l}+\sqrt{r_jr_l}=\sqrt{r_jr_k}$, we have
$$r_i=-\frac{C}{B}>0.$$
This corresponds to the case that there is a straight line tangent to the circles $C_j, C_k, C_l$ on the same side.
See (b) in Figure \ref{circle}.

  \item[(3)] In the other cases, we have $A<0$, $B>0$ and $C<0$ with discriminant $\Delta>0$.
Then the quadratic equation (\ref{quadratic equation in r_i}) have two different positive roots
$\frac{-B+\sqrt{\Delta}}{2A}$ and $\frac{-B-\sqrt{\Delta}}{2A}$.
This corresponds to the case that there are two circles
with finite radii externally tangent to $C_j, C_k, C_l$ simultaneously.
See (c) in Figure \ref{circle}.
Without loss of generality, we assume $r_j\leq r_k\leq r_l$. We claim that $\frac{-B-\sqrt{\Delta}}{2A}>r_j$.
Then $r_i=\frac{-B+\sqrt{\Delta}}{2A}$ by $r_i<\min\{r_j, r_k, r_l\}$.

To prove the claim, note that
\begin{equation*}
\begin{aligned}
A=r_j^2r_k^2r_l^2&(\frac{1}{\sqrt{r_j}}+\frac{1}{\sqrt{r_k}}+\frac{1}{\sqrt{r_l}})(\frac{1}{\sqrt{r_j}}+\frac{1}{\sqrt{r_k}}-\frac{1}{\sqrt{r_l}})\\
&(\frac{1}{\sqrt{r_j}}+\frac{1}{\sqrt{r_l}}-\frac{1}{\sqrt{r_k}})(\frac{1}{\sqrt{r_k}}+\frac{1}{\sqrt{r_l}}-\frac{1}{\sqrt{r_j}}).
\end{aligned}
\end{equation*}
$A<0$ implies $\frac{1}{\sqrt{r_j}}>\frac{1}{\sqrt{r_k}}+\frac{1}{\sqrt{r_l}}$.
To simplify the notations, set $a=\frac{1}{\sqrt{r_j}}$, $b=\frac{1}{\sqrt{r_k}}$, $c=\frac{1}{\sqrt{r_l}}$, then $a>b+c$.

$\frac{-B-\sqrt{\Delta}}{2A}>r_j$ if and only if
\begin{equation*}
\begin{aligned}
2 r_jr_kr_l&(r_kr_l+r_jr_l+r_jr_k)+\sqrt{16r_j^3r_k^3r_l^3(r_j+r_k+r_l)}\\
&>-2r_j(2r_jr_kr^2_l+2r_jr_lr^2_k+2r_kr_lr^2_j-r_k^2r_l^2-r_j^2r_l^2-r_j^2r_k^2),
\end{aligned}
\end{equation*}
if and only if
\begin{equation*}
\begin{aligned}
 \frac{1}{r_j}+\frac{1}{r_k}+&\frac{1}{r_l}+2\sqrt{\frac{1}{r_jr_k}+\frac{1}{r_jr_l}+\frac{1}{r_kr_l}}\\
&>r_j(\frac{1}{r_j^2}+\frac{1}{r_k^2}+\frac{1}{r_l^2}-\frac{2}{r_jr_k}-\frac{2}{r_jr_l}-\frac{2}{r_kr_l}),
\end{aligned}
\end{equation*}
if and only if
\begin{equation*}
\begin{aligned}
a^2+b^2+&c^2+2\sqrt{a^2b^2+a^2c^2+b^2c^2}\\
&>\frac{1}{a^2}(a^4+b^4+c^4-2a^2b^2-2a^2c^2-2b^2c^2),
\end{aligned}
\end{equation*}
if and only if
\begin{equation}\label{inequality}
\begin{aligned}
3a^2(b^2+c^2)+2a^2\sqrt{a^2b^2+a^2c^2+b^2c^2}>(b^2-c^2)^2.
\end{aligned}
\end{equation}
By $a>b+c$, we have
$$3a^2(b^2+c^2)>3(b+c)^2(b^2+c^2)>(b+c)^2(b-c)^2=(b^2-c^2)^2,$$
which implies (\ref{inequality}).
This completes the proof of the claim and the lemma.\hfill\rule{1.5mm}{3mm}
\end{description}

Note that $\partial \Omega_{jkl}=\overline{\Omega}_{jkl}\setminus \Omega_{jkl}=\{(r_j,r_k,r_l)\in \mathbb{R}^3_{>0}|A=0\}$.
It is straightforward to show that, as $(r_j,r_k,r_l)$ tends to a point in $\partial \Omega_{jkl}$,
$$\frac{-B+\sqrt{B^2-4AC}}{2A}\rightarrow -\frac{C}{B},$$
which implies that $r_i=f(r_j, r_k, r_l)$ in (\ref{equation of r_i}) is a continuous and
piecewise analytic function of $(r_j, r_k, r_l)\in \mathbb{R}^3_{>0}$.

By (\ref{derivative of Q}), the degenerate set $V_i$ is
$$V_i=\{(r_i,r_j,r_k,r_l)\in \mathbb{R}^4_{>0}|0<r_i\leq f(r_j,r_k,r_l)\},$$
which is a simply connected subset of $\mathbb{R}^4_{>0}$.
Similarly, we have
\begin{equation*}
\begin{aligned}
V_j=&\{(r_i,r_j,r_k,r_l)\in \mathbb{R}^4_{>0}|0<r_j\leq f(r_i,r_k,r_l)\},\\
V_k=&\{(r_i,r_j,r_k,r_l)\in \mathbb{R}^4_{>0}|0<r_k\leq f(r_i,r_j,r_l)\},\\
V_l=&\{(r_i,r_j,r_k,r_l)\in \mathbb{R}^4_{>0}|0<r_l\leq f(r_i,r_j,r_k)\}.
\end{aligned}
\end{equation*}
Therefore, we have
$$\mathbb{R}^4_{>0}=\Omega_{ijkl}\cup V_i\cup V_j \cup V_k \cup V_l.$$

We claim that $V_i, V_j, V_k, V_l$ are mutually disjoint.
Otherwise that $V_i\cap V_j\neq \emptyset$
and $r=(r_i, r_j, r_k, r_l)\in V_i\cap V_j \subset \mathbb{R}^4_{>0}$.
By the geometric meaning of the critical degenerate case, we have
\begin{equation}\label{i in jkl}
A_{\triangle ijk}+A_{\triangle ijl}+A_{\triangle ikl}\leq A_{\triangle jkl}
\end{equation}
and
\begin{equation}\label{j in ikl}
A_{\triangle ijk}+A_{\triangle ijl}+A_{\triangle jkl}\leq A_{\triangle ikl},
\end{equation}
where $A_{\triangle ijk}$ denotes the area of the triangle $\{ijk\}\in F$ with edge lengths $l_{ij}=r_i+r_j, l_{ik}=r_i+r_k, l_{jk}=r_j+r_k$.
Combining (\ref{i in jkl}) with (\ref{j in ikl}), we have $A_{\triangle ijk}+A_{\triangle ijl}\leq 0$,
which is impossible.
So we have $V_i\cap V_j=\emptyset$.
This completes the proof for the theorem with Euclidean background geometry.

The proof for the case of hyperbolic background geometry is similar.
The boundary of $V_i$ in $\mathbb{R}^4_{>0}$ is given by the function
\begin{equation*}
\begin{aligned}
\tanh r_i=\left\{
            \begin{array}{ll}
              \frac{-\mathcal{B}+\sqrt{\mathcal{B}^2-4\mathcal{A}\mathcal{C}}}{2\mathcal{A}}, & \hbox{$(r_j, r_k, r_l)\in\mathbb{R}^3_{>0}\setminus\partial \widetilde{\Omega}_{jkl}$,} \\
              -\frac{\mathcal{C}}{\mathcal{B}}, & \hbox{$(r_j, r_k, r_l)\in \partial \widetilde{\Omega}_{jkl}$,}
            \end{array}
          \right.
\end{aligned}
\end{equation*}
where
\begin{equation*}
\begin{aligned}
\mathcal{A}=&4\tanh^2 r_j \tanh^2 r_k\tanh^2 r_l\\
&+\left(\tanh r_j \tanh r_k+\tanh r_j\tanh r_l+\tanh r_k\tanh r_l\right)^2\\
&-2\left(\tanh^2 r_j \tanh^2 r_k+\tanh^2 r_j \tanh^2 r_l+\tanh^2 r_k \tanh^2 r_l\right),\\
\mathcal{B}=& 2\tanh r_j \tanh r_k\tanh r_l(\tanh r_j \tanh r_k+\tanh r_j \tanh r_l+\tanh r_k \tanh r_l),\\
\mathcal{C}=&-\tanh^2 r_j \tanh^2 r_k \tanh^2 r_l,\\
\partial \widetilde{\Omega}_{jkl}=&\{(r_j, r_k, r_l)\in \mathbb{R}^3_{>0}|\mathcal{A}=0\}.
\end{aligned}
\end{equation*}
Set $r_i=g(r_j, r_k, r_l)$, then $g$ is continuous and piecewise analytic. The corresponding degenerate set $V_i$ is
$$V_i=\{(r_i, r_j, r_k, r_l)\in \mathbb{R}^4_{>0}|0<r_i\leq g(r_j, r_k, r_l)\}.$$
The rest of the proof for the hyperbolic case is similar to that of the Euclidean case, so we omit the details here.
\qed

\begin{remark}
By Theorem \ref{structure of admissible r}, the admissible space $\Omega_{ijkl}$ of sphere
packing metrics for a single tetrahedron is homotopy equivalent to $\mathbb{R}^4_{>0}$.
This provides another proof of Cooper-Rivin's Lemma \ref{simply connectness of ijkl} \cite{CR}
that $\Omega_{ijkl}$ is simply connected.
\end{remark}

In the following, we take $\partial_i \Omega_{ijkl}=\overline{\Omega}_{ijkl}\cap V_i$.

\section{Cooper-Rivin's action functional and its extension}\label{Section 3}
\subsection{Cooper-Rivin's action functional}

For a triangulated 3-manifold $(M, \mathcal{T})$ with sphere packing metric $r$,
Cooper and Rivin \cite{CR} introduced the definition (\ref{CR curvature}) of combinatorial scalar curvature $K_{i}$
at the vertex $i$
$$K_{i}= 4\pi-\sum_{\{i,j,k,l\}\in T}\alpha_{ijkl},$$
where $\alpha_{ijkl}$ is the solid angle at the vertex $i$ of the tetrahedron $\{i,j,k,l\}\in T$
and the sum is taken over all tetrahedra with $i$ as one of its vertices.
Given a single tetrahedron $\{ijkl\}\in T$, we usually denote the solid angle $\alpha_{ijkl}$ at the
vertex $v_i$ by $\alpha_i$ for simplicity.
$K_i$ locally measures the difference between the volume growth rate of a small ball
centered at vertex $v_i$ in $M$ and a Euclidean ball of the same radius.
Cooper and Rivin's definition (\ref{CR curvature}) of combinatorial scalar curvature is motivated by the
fact that, in the smooth case, the scalar curvature at a point $P$ locally measures the difference of
the volume growth rate of the geodesic ball with center $P$ to the Euclidean ball with the same radius \cite{Be, LP}.
In fact, for a geodesic ball $B(P, r)$ in an n-dimensional Riemannian manifold $(M^n, g)$ with center $P$ and radius $r$, we have
the following asymptotical expansion for the volume of $B(P, r)$
$$\text{Vol}(B(P, r))=\omega(n)r^n\left(1-\frac{1}{6(n+2)}R(P)r^2+o(r^2)\right),$$
where $\omega(n)$ is the volume of the unit ball in $\mathbb{R}^n$ and $R(P)$ is the scalar curvature of $(M, g)$ at $P$.
From this point of view, Cooper and Rivin's definition of combinatorial scalar curvature is a
good candidate for combinatorial scalar curvature with geometric meaning similar to the smooth case.

\begin{lemma}[\cite{CR, G1, G2, R}]\label{concave for one tetrahedron}
For a tetrahedron $\{ijkl\}\in T$, set
\begin{equation*}
\begin{aligned}
\mathcal{S}_{ijkl}=\left\{
                     \begin{array}{ll}
                       \sum_{\mu\in \{i,j,k,l\}}\alpha_{\mu}r_{\mu}, & \hbox{Euclidean background geometry} \\
                       \sum_{\mu\in \{i,j,k,l\}}\alpha_{\mu}r_{\mu}+2 \text{Vol}, & \hbox{hyperbolic background geometry}
                     \end{array}
                   \right.
\end{aligned},
\end{equation*}
where Vol denotes the volume of the tetrahedron for the hyperbolic background geometry.
Then
\begin{equation}\label{schlafli formula}
d\mathcal{S}_{ijkl}=\sum_{\mu\in \{i,j,k,l\}}\alpha_{\mu}dr_{\mu}=\alpha_idr_i+\alpha_jdr_j+\alpha_kdr_k+\alpha_ldr_l.
\end{equation}
Furthermore, the Hessian of $\mathcal{S}_{ijkl}$ is negative semi-definite with kernel
$\{t(r_i, r_j, r_k, r_l)|t\in \mathbb{R}\}$ for the Euclidean background geometry and negative definite for the hyperbolic background geometry.
\end{lemma}

(\ref{schlafli formula}) implies that $\alpha_idr_i+\alpha_jdr_j+\alpha_kdr_k+\alpha_ldr_l$
is a closed $1$-form on $\Omega_{ijkl}$.
Combining Lemma \ref{simply connectness of ijkl} with Lemma \ref{concave for one tetrahedron},
we have

\begin{lemma}[\cite{CR, G1,  G2, R}]\label{F_ijkl}
Given a Euclidean or hyperbolic tetrahedron $\{ijkl\}\in T$ and $r_0\in \Omega_{ijkl}$,
\begin{equation}\label{definition of F_ijkl}
F_{ijkl}(r)=\int_{r_0}^r\alpha_idr_i+\alpha_jdr_j+\alpha_kdr_k+\alpha_ldr_l
\end{equation}
is a well-defined locally concave function on $\Omega_{ijkl}$.
Furthermore, $F_{ijkl}(r)$ is strictly concave on $\Omega_{ijkl}\cap\{r_i^2+r_j^2+r_k^2+r_l^2=c\}$ for any $c>0$ in the Euclidean background geometry
and strictly concave on $\Omega_{ijkl}$ in the hyperbolic background geometry.
\end{lemma}

\begin{remark}
It is observed \cite{CR, G1} that (\ref{schlafli formula}) is essentially the Schl\"{a}fli formula.
\end{remark}

Using Lemma \ref{concave for one tetrahedron}, we have the following property for the combinatorial scalar curvature.

\begin{lemma}\label{property of Lambda}
(\cite{CR, G1, G2, R})
Suppose $(M, \mathcal{T})$ is a triangulated 3-manifold with sphere packing metric $r$,
$\mathcal{S}$ is the total combinatorial scalar curvature defined as
\begin{equation*}
\begin{aligned}
\mathcal{S}(r)=\left\{
                 \begin{array}{ll}
                   \sum K_ir_i, & \hbox{Euclidean background geometry} \\
                   \sum K_ir_i-2\text{Vol}(M), & \hbox{hyperbolic background geometry}
                 \end{array}.
               \right.
\end{aligned}
\end{equation*}
Then we have
\begin{equation*}
d\mathcal{S}=\sum_{i=1}^NK_idr_i.
\end{equation*}
Set
\begin{equation}\label{Matrix Lambda}
\begin{aligned}
\Lambda=\operatorname{Hess}_r\mathcal{S}=
\frac{\partial(K_{1},\cdots,K_{N})}{\partial(r_{1},\cdots,r_{N})}=
\left(
\begin{array}{ccccc}
 {\frac{\partial K_1}{\partial r_1}}& \cdot & \cdot & \cdot &  {\frac{\partial K_1}{\partial r_N}} \\
 \cdot & \cdot & \cdot & \cdot & \cdot \\
 \cdot & \cdot & \cdot & \cdot & \cdot \\
 \cdot & \cdot & \cdot & \cdot & \cdot \\
 {\frac{\partial K_N}{\partial r_1}}& \cdot & \cdot & \cdot &  {\frac{\partial K_N}{\partial r_N}}
\end{array}
\right).
\end{aligned}
\end{equation}
In the case of Euclidean background geometry, $\Lambda$ is symmetric and positive semi-definite with rank $N-1$ and
kernel $\{tr|t\in\mathbb{R}\}$.
In the case of hyperbolic background geometry, $\Lambda$ is symmetric and positive definite.
\end{lemma}

We refer the readers to \cite{G1, G4} for a nice geometrical explanation of $\frac{\partial K_i}{\partial r_j}$.
It should be emphasized that, as pointed out by Glickenstein \cite{G2}, the elements $\frac{\partial K_i}{\partial r_j}$
for $i\sim j$ may be negative, which is different from two-dimensional case.

\subsection{Extension of Cooper-Rivin's action functional}
For a single Euclidean or hyperbolic tetrahedron $\{ijkl\}\in T$, the solid angle function
$\alpha(r)=(\alpha_i(r), \alpha_j(r), \alpha_k(r), \alpha_l(r))$ is defined on the admissible space $\Omega_{ijkl}$.
We will extend the solid angle function
to $\mathbb{R}^4_{>0}$ continuously by making it to be a constant function
on each connected components of $\mathbb{R}^4_{>0}\setminus\Omega_{ijkl}$,
which is called a continuous extension by constants in \cite{L3}.
We have the following result.
\begin{lemma}\label{extension of solid angle}
The solid angle function $\alpha=(\alpha_i, \alpha_j, \alpha_k, \alpha_l)$ defined on $\Omega_{ijkl}$
can be extended continuously by constants to a function
$\widetilde{\alpha}=(\widetilde{\alpha}_i, \widetilde{\alpha}_j, \widetilde{\alpha}_k, \widetilde{\alpha}_l)$ defined on $\mathbb{R}^4_{>0}$.
\end{lemma}

\proof
The extension $\widetilde{\alpha}_i$ of $\alpha_i$ is defined to be $\widetilde{\alpha}_i(r)=2\pi$ for $r=(r_i, r_j, r_k, r_l)\in V_i$ and
$\widetilde{\alpha}_i(r)=0$ for $r=(r_i, r_j, r_k, r_l)\in V_\alpha$ with $\alpha\in \{j, k, l\}$.
The extensions $\widetilde{\alpha}_j$,
$\widetilde{\alpha}_k$, $\widetilde{\alpha}_l$ of $\alpha_j$, $\alpha_k$, $\alpha_l$ are defined similarly.

If $r=(r_i, r_j, r_k, r_l)\in \Omega_{ijkl}$ and $r \rightarrow P$ for some point $P\in \partial_i\Omega_{ijkl}$,
then geometrically the tetrahedron $\{ijkl\}$ tends
to degenerate with the center $v_i$ of the sphere $S_i$ tends to lie in the geodesic plane defined by $v_j, v_k, v_l$ and
the corresponding circle $C_i$ is externally tangent to $C_j, C_k, C_l$ with $i$ in the interior of the triangle $\triangle jkl$.
Then by Proposition \ref{degenerate lemma}, we have
$\alpha_i\rightarrow 2\pi, \alpha_j\rightarrow 0, \alpha_k\rightarrow 0$ and $\alpha_l\rightarrow 0$, as $r\rightarrow P\in \partial_i\Omega_{ijkl}$.
This implies that the extension $\widetilde{\alpha}$ of $\alpha$ is continuous on $\mathbb{R}^4_{>0}$.\qed

Before going on, we recall the following definition and theorem of Luo in \cite{L3}.
\begin{definition}
A differential 1-form $w=\sum_{i=1}^n a_i(x)dx^i$ in an open set $U\subset \mathbb{R}^n$ is said to be continuous if
each $a_i(x)$ is continuous on $U$.  A continuous differential 1-form $w$ is called closed if $\int_{\partial \tau}w=0$ for each
triangle $\tau\subset U$.
\end{definition}

\begin{theorem}[\cite{L3}, Corollary 2.6]\label{Luo's convex extention}
Suppose $X\subset \mathbb{R}^n$ is an open convex set and $A\subset X$ is an open subset of $X$ bounded by a $C^1$
smooth codimension-1 submanifold in $X$. If $w=\sum_{i=1}^na_i(x)dx_i$ is a continuous closed 1-form on $A$ so that
$F(x)=\int_a^x w$ is locally convex on $A$ and each $a_i$ can be extended continuous to $X$ by constant functions to a
function $\widetilde{a}_i$ on $X$, then  $\widetilde{F}(x)=\int_a^x\sum_{i=1}^n\widetilde{a}_i(x)dx_i$ is a $C^1$-smooth
convex function on $X$ extending $F$.
\end{theorem}

Combining Theorem \ref{structure of admissible r}, Lemma \ref{F_ijkl}, Lemma \ref{extension of solid angle} and Theorem \ref{Luo's convex extention},
we have

\begin{lemma}\label{extension of F_ijkl}
For a Euclidean or hyperbolic tetrahedron $\{i,j,k,l\}\in T$,
the function $F_{ijkl}(r)$ defined on $\Omega_{ijkl}$ in (\ref{definition of F_ijkl}) can be extended to
\begin{equation}\label{extended Ricci potential for a tetrahedron}
\widetilde{F}_{ijkl}(r)=\int_{r_0}^r\widetilde{\alpha}_idr_i+\widetilde{\alpha}_jdr_j+\widetilde{\alpha}_kdr_k+\widetilde{\alpha}_ldr_l,
\end{equation}
which is a $C^1$-smooth concave function defined on $\mathbb{R}^4_{>0}$ with
$$\nabla_r\widetilde{F}_{ijkl}=\widetilde{\alpha}^T=(\widetilde{\alpha}_i, \widetilde{\alpha}_j, \widetilde{\alpha}_k, \widetilde{\alpha}_l)^T.$$
\end{lemma}

\section{Proof of the global rigidity}\label{Section 4}
\subsection{Proof of Theorem \ref{main theorem rigidity}}
We first introduce the following definition.
\begin{definition}
Given a function $f\in \mathbb{R}^V$, if there exists a sphere packing metric $r\in \Omega$ such that
$K(r)=f$, then $f$ is called an admissible curvature function.
\end{definition}

Theorem \ref{main theorem rigidity} is equivalent to the following form.

\begin{theorem}
Suppose $\overline{K}$ is an admissible curvature function
on a connected closed triangulated 3-manifold $(M, \mathcal{T})$,
then there exists only one admissible sphere packing metric in $\Omega$ with combinatorial scalar curvature $\overline{K}$
(up to scaling in the case of Euclidean background geometry).
\end{theorem}
\proof
We only prove the Euclidean sphere packing case and the hyperbolic sphere packing case is proved similarly.
Suppose $r_0\in \Omega$ is a sphere packing metric.
Define a Ricci potential function
\begin{equation}\label{Ricci potential}
\begin{aligned}
\widetilde{F}(r)=-\sum_{\{ijkl\}\in T}\widetilde{F}_{ijkl}(r_i, r_j, r_k, r_l)+\sum_{i=1}^N(4\pi-\overline{K}_i)r_i,
\end{aligned}
\end{equation}
where the function $\widetilde{F}_{ijkl}$ is defined by (\ref{extended Ricci potential for a tetrahedron}).
Note that the second term in the right-hand-side of (\ref{Ricci potential}) is linear in $r$ and well-defined on $\mathbb{R}^N_{>0}$.
Combining with Lemma \ref{extension of F_ijkl},
we have $\widetilde{F}(r)$ is a well-defined $C^1$-smooth convex function on $\mathbb{R}^N_{>0}$.
Furthermore,
\begin{equation}\label{gradient of Ricci potential}
\nabla_{r_i}\widetilde{F}=-\sum_{\{ijkl\}\in T}\widetilde{\alpha}_{ijkl}+(4\pi-\overline{K}_i)=\widetilde{K}_i-\overline{K}_i,
\end{equation}
where $\widetilde{K}_i=4\pi-\sum_{\{ijkl\}\in T}\widetilde{\alpha}_{ijkl}$ is an extension of $K_i=4\pi-\sum_{\{ijkl\}\in T}\alpha_{ijkl}$.

If there are two different sphere packing metrics $\overline{r}_A$ and $\overline{r}_B$
in $\Omega$ with the same combinatorial scalar curvature $\overline{K}$,
then
$$\nabla \widetilde{F}(\overline{r}_A)=\nabla \widetilde{F}(\overline{r}_B)=0$$
by (\ref{gradient of Ricci potential}).
Set
\begin{equation*}
\begin{aligned}
f(t)=&\widetilde{F}((1-t)\overline{r}_A+t\overline{r}_B)\\
=&\sum_{\{ijkl\}\in T}f_{ijkl}(t)+\sum_{i=1}^N(4\pi-\overline{K}_i)[(1-t)\overline{r}_{A, i}+t\overline{r}_{B, i}],
\end{aligned}
\end{equation*}
where
\begin{equation*}
\begin{aligned}f_{ijkl}(t)=-\widetilde{F}_{ijkl}(&(1-t)\overline{r}_{A, i}+t\overline{r}_{B, i}, (1-t)\overline{r}_{A, j}+t\overline{r}_{B, j}, \\ &(1-t)\overline{r}_{A, k}+t\overline{r}_{B, k}, (1-t)\overline{r}_{A, l}+t\overline{r}_{B, l}).
\end{aligned}
\end{equation*}
Then $f(t)$ is a $C^1$-smooth convex function on $[0, 1]$ and $f'(0)=f'(1)=0$. Therefore
$f'(t)\equiv 0$ on $[0, 1]$.

Note that $\overline{r}_A\in \Omega$ and $\Omega$ is an open subset of $\mathbb{R}^N_{>0}$,
there exists $\epsilon>0$ such that $(1-t)\overline{r}_{A}+t\overline{r}_{B}\in \Omega$ for $t\in [0, \epsilon]$.
Hence $f(t)$ is $C^2$ (in fact $C^{\infty}$) on $[0, \epsilon]$.
$f'(t)\equiv 0$ on $[0, 1]$ implies that $f''(t)\equiv 0$ on $[0, \epsilon]$. But for $t\in [0, \epsilon]$, we have
$$f''(t)
=(\overline{r}_A-\overline{r}_B)\Lambda(\overline{r}_A-\overline{r}_B)^T,$$
where $\Lambda$ is the matrix defined in (\ref{Matrix Lambda}).
By Lemma \ref{property of Lambda}, we have $\overline{r}_A=c\overline{r}_B$ for some positive constant $c\in \mathbb{R}$.
So there exists only one Euclidean sphere packing metric up to scaling with combinatorial scalar curvature equal to $\overline{K}$.   \qed

\begin{remark}
The proof of Theorem \ref{main theorem rigidity} is based on a variational principle introduce by Colin de Verdi\`{e}re \cite{DV}.
Bobenko, Pinkall and Springborn \cite{BPS} introduced a method to extend a locally convex function on a nonconvex domain to a convex function
and solved affirmably a conjecture of Luo \cite{L1} on the global rigidity of piecewise linear metrics on surfaces.
Using the method of extension,
Luo \cite{L3} proved the global rigidity of inversive distance circle packing metrics for nonnegative inversive distance
and the author \cite{X} proved the global rigidity of inversive distance circle packing metrics when the inversive distance is in $(-1, +\infty)$.
The method of extension was further used to study the deformation of combinatorial curvatures on surfaces in \cite{GJ0, GJ1, GJ2, GJ3, GX6}.
\end{remark}

\subsection{Proof of Theorem \ref{main theorem global rigidity for alpha curvature}}
The proof is the same as that of Theorem \ref{main theorem rigidity} using a similar defined Ricci potential function.
To be more precisely, for the Euclidean sphere packing metrics, define
\begin{equation*}
\begin{aligned}
F(r)=-\sum_{\{ ijkl\}\in T}F_{ijkl}(r_i, r_j, r_k, r_l)+\int_{r_0}^r\sum_{i=1}^N(4\pi-\overline{R}_ir_i^{\alpha})dr_i
\end{aligned}
\end{equation*}
on $\Omega$, where $r_0\in \Omega$ and $F_{ijkl}$ is defined by (\ref{definition of F_ijkl}).
$F(r)$ can be extended to a $C^1$-smooth function
$$\widetilde{F}(r)=-\sum_{\{ ijkl\}\in T}\widetilde{F}_{ijkl}(r_i, r_j, r_k, r_l)+\int_{r_0}^r\sum_{i=1}^N(4\pi-\overline{R}_ir_i^\alpha)dr_i$$
defined on $\mathbb{R}^N_{>0}$ by Lemma \ref{extension of F_ijkl},
where $\widetilde{F}_{ijkl}$ is defined by (\ref{extended Ricci potential for a tetrahedron}).

Following the same arguments as that in the proof of Theorem \ref{main theorem rigidity}, there exists $\epsilon\in (0,1)$ such that
\begin{equation}\label{second derivative of f_alpha}
\begin{aligned}
f''(t)=(\overline{r}_A-\overline{r}_B)\cdot \operatorname{Hess}_r F \cdot (\overline{r}_A-\overline{r}_B)^T\equiv 0
\end{aligned}
\end{equation}
for $t\in[0, \epsilon]$, where
\begin{equation*}
\begin{aligned}
\operatorname{Hess}_r F=\Lambda-\alpha\left(
                   \begin{array}{ccc}
\overline{R}_1r_1^{\alpha-1} &   &   \\
                       & \ddots &   \\
                       &   & \overline{R}_Nr_N^{\alpha-1}\\
                   \end{array}
                 \right).
\end{aligned}
\end{equation*}

In the case that $\alpha\overline{R}\equiv0$,
$f''(t)=(\overline{r}_A-\overline{r}_B) \Lambda  (\overline{r}_A-\overline{r}_B)^T=0$
for $t\in [0, \epsilon]$. By Lemma \ref{property of Lambda}, we have $\overline{r}_A=c\overline{r}_B$ for some positive constant $c\in \mathbb{R}$.
So there exists at most one Euclidean sphere packing metric
with combinatorial $\alpha$-curvature equal to $\overline{R}$ up to scaling.

In the case that $\alpha\overline{R}\leq 0$ and $\alpha\overline{R}\not\equiv0$,
$\operatorname{Hess}_r F$ is positive definite on $\Omega$ by Lemma \ref{property of Lambda}.
Then (\ref{second derivative of f_alpha}) implies $\overline{r}_A=\overline{r}_B$.
Therefore there exists at most one Euclidean sphere packing
metric with combinatorial $\alpha$-curvature equal to $\overline{R}$.

For the hyperbolic case,
\begin{equation*}
\begin{aligned}
F(r)=-\sum_{\{ ijkl\}\in T}F_{ijkl}(r_i, r_j, r_k, r_l)+\int_{r_0}^r\sum_{i=1}^N(4\pi-\overline{R}_i\tanh^\alpha\frac{r_i}{2})dr_i.
\end{aligned}
\end{equation*}
The proof is similar to the Euclidean case,
so we omit the details here.\qed
~\\

\textbf{Acknowledgements}

~\\
The research of the author is supported by National Natural Science Foundation of China under
grant no. 11301402, Fundamental Research Funds for the Central Universities under grant no. 2042018kf0246
and Hubei Provincial Natural Science Foundation of China under grant no. 2017CFB681.
The author would like to thank Professor Feng Luo for his interest in this work and careful reading of the paper
and Professor Kehua Su for the help of Figure \ref{circle}.
The author also would like to thank the referees for their valuable suggestions that greatly improve the paper.

(Xu Xu)

$^1$ School of Mathematics and Statistics, Wuhan University, Wuhan 430072, P.R. China

$^2$ Hubei Key Laboratory of Computational Science, Wuhan University, Wuhan 430072, P.R.China

E-mail: xuxu2@whu.edu.cn\\[2pt]

\end{document}